\documentclass[12pt]{article}
\usepackage{amssymb}
\newtheorem{theorem}{\sc Theorem}[section]
\newtheorem{definition}{\sc Definition}[section]
\newtheorem{lemma}{\sc Lemma}[section]
\newtheorem{corollary}{\sc Corollary}[section]
\newtheorem{example}{\sc Example}[section]

\begin{document}
\baselineskip=24 pt
\title{\bf Asymptotic equivalence of differential equations and asymptotically almost periodic solutions}
\author{M. U. Akhmet$^1$\thanks{Corresponding author. M.U. Akhmet is previously known as M. U. Akhmetov.},  M.A. Tleubergenova$^2$
and  A. Zafer$^3$ }
\date{{\small$^1$ Department of Mathematics and Institute of Applied Mathematics, Middle East
Technical University, 06531 Ankara, Turkey, marat@metu.edu.tr\\$^2$ 
Department of Mathematics, K. Zhubanov Aktobe State Pedagogical University, 463000 Aktobe, pr. Moldagulovoy, 34, Kazakhstan, \\
$^3$ Department of Mathematics, Middle East
Technical University, 06531 Ankara, Turkey, zafer@metu.edu.tr }}

\maketitle



\newpage
\noindent {\bf Proposed running head:} Asymptotic equivalence, almost periodic solutions\\

\begin{abstract}
In this  paper we establish  asymptotic (biasymptotic) equivalence
between  spaces of solutions of a given linear homogeneous system  and a perturbed system.
The perturbations are of either linear or weakly linear characters.
Existence of a homeomorphism between subspaces of
almost periodic and  asymptotically (biasymptotically) almost periodic solutions is also obtained.
\end{abstract}
\maketitle
\section{Introduction and Preliminaries}
Let $\mathbb N$ and $\mathbb R$ be  sets of all natural and real
numbers, respectively. Denote by $||\cdot||$ the Euclidean norm in
$\mathbb R^n$, $n \in \mathbb N$, and by $C(\mathbb X,\mathbb Y)$
the space of all continuous functions defined on $\mathbb X$ with
values in $\mathbb Y$.

We shall recall the definitions of almost periodic and
asymptotically almost periodic functions, see \cite{amerio, cord, fink} for more
details.

A number  $\tau \in \mathbb R$ is called an $\epsilon -$ translation number of a function
$f\in C(\mathbb R,\mathbb R^n)$  if  $|| f(t+\tau) - f(t)||< \epsilon$ for all $t \in \mathbb R.$
A function $f \in C (\mathbb R,\mathbb R^n)$ is called almost periodic
if for a given  $\epsilon \in \mathbb R$,  $\epsilon > 0$, there exists  a relatively dense set of
$\epsilon -$ translation numbers of $f$. The set of all almost periodic  functions is denoted by
${\cal AP}(\mathbb R)$.

A function $f \in C (\mathbb R,\mathbb R^n)$ is called asymptotically almost periodic
if there is a function $g \in {\cal AP}(R)$ and a function  $\phi \in C(\mathbb R,\mathbb R^n)$ with
$\lim_{t \rightarrow \infty} \phi(t) = 0$ such that
  $f(t) = g(t) + \phi(t).$

The basic definition of an almost periodic function given by H. Bohr has been modified by several
authors \cite{amerio, cord, fink, hw, sp, z}.  Below we  introduce a new notion with regard to almost
periodic functions.

\begin{definition} A function $f \in C (\mathbb R,\mathbb R^n)$ is called biasymptotically almost periodic
  if $f(t) = g(t) + \phi(t)$ for some $g \in {\cal AP}(\mathbb R)$ and $\phi\in C(\mathbb R,\mathbb R^n)$ with $\lim_{t \rightarrow \pm \infty} \phi(t) = 0.$
 \end{definition}
Note that every biasymptotically almost periodic function  is a
pseudo almost periodic function \cite{z}, but not conversely.

In this paper we are concerned with the linear system
\begin{eqnarray}
&& y' = [A(t) + B(t)]y,
\label{1}
\end{eqnarray}
which may be viewed  as a perturbation of
\begin{eqnarray}
&& x' = A(t)x,
\label{2}
\end{eqnarray}
where $x,y \in \mathbb R^n$, and $A, B \in C(\mathbb R,\mathbb R^{n\times n})$.

Moreover we consider the quasilinear systems of the form
\begin{eqnarray}
&& y' = C y + f(t,y)
\label{e2n}
\end{eqnarray}
and the corresponding homogeneous linear system
\begin{eqnarray}
&& x' = C x,
\label{e1n}
\end{eqnarray}
where $x,y \in  \mathrm R^n$, $C\in \mathbb R^{n\times n}$, and $f(t,x) \in C( \mathbb R \times \mathbb R^n,\mathbb R^n)$ such that $$f(t,0)=0\quad \mbox{for all $t\in \mathbb R$.}$$
\begin{definition}[\cite{c,ns}]
A  homeomorphism  between solutions $x(t)$  and   $y(t)$
is called an asymptotic equivalence  if   $x(t) - y(t) \rightarrow 0$ as $ t\rightarrow  \infty.$
\label{defn1}
\end{definition}
\begin{definition}
A  homeomorphism between solutions $x(t)$  and   $y(t)$
is called an biasymptotic equivalence  if   $x(t) - y(t) \rightarrow 0$ as $ t\rightarrow \pm \infty.$
\label{defn1prime}
\end{definition}

Our main objective is to investigate the problem of  asymptotic  equivalence of
systems and to prove the existence of asymptotically and biasymptotically almost periodic solutions
of (\ref{1}) and  (\ref{e2n}) .

The classical theorem of Levinson \cite{nl} states that if the trivial solution of (\ref{2}) is uniformly stable,
 $A(t)\equiv A$, and
\begin{equation}\int_{0}^{\infty}\mid\mid B(t)\mid\mid\,dt <  \infty, \label{b1}\end{equation}
then (\ref{1}) and (\ref{2}) are asymptotically equivalent.  In
the case when $A$ is not a constant matrix,  Wintner \cite{aw}
proved that the  above conclusion  remains valid if all solutions
of (\ref{2}) are bounded, (\ref{b1}) is satisfied, and
$$\liminf_{t\to\infty}\int_{0}^{t}\mbox{Trace}[B(s)]\,ds > - \infty.$$
Later, Yakubovi\v c \cite{ya} considered (\ref{e2n}) and obtained asymptotic equivalence of (\ref{e1n}) and (\ref{e2n}),
see \cite{ns, ya} for details.
After the pioneering works of Levinson, Wintner, and Yakubovi\v c,  the problem
of asymptotic equivalence of  differential systems including linear, nonlinear, and functional equations has
been investigated by many authors; see e.g. [2, 12--20],
and the references cited therein.
Two interesting articles in this direction which also motivate our study here in this paper were written by  M. R\'ab \cite{mr1,mr2}.
In fact, the main result in \cite{mr2} is an improvement of the earlier one in \cite{mr1}, which
we have employed in  our work.

Asymptotically almost periodic functions were introduced by
Fr\'echet \cite{f,f1}.   The  existence of  this type of solutions  was investigated by 
A.M. Fink  \cite{fink} (Theorem 9.5)  for the first time. 
For more results on the existence of asymptotically almost
periodic solutions of different types of equations we refer to
\cite{ h1,h2,mt, pc, v, y} and the references cited therein. 
In this work we exploit the idea of A.M. Fink  
to obtain the existence of asymptotically almost periodic solutions of linear and quasilinear 
systems. Moreover, we prove a theorem about
biasymptotic equivalence of linear systems and a theorem on the
existence of biasymptotically almost periodic solutions.
Apparently the notions of biasymptotic equivalence and a
biasymptotically almost periodic function are introduced for the
first time in this paper.

The paper is organized as follows. In the next section, we prove a main lemma of  R\'ab and obtain sufficient conditions concerning
the asymptotic equivalence of (\ref{1}) and (\ref{2}), and the existence of a family of asymptotically almost periodic
solutions of the system   (\ref{1}). The third section  is devoted to the problem of the asymptotic equivalence of
systems  (\ref{e1n}) and ({\ref{e2n}) and the problem of existence of asymptotically almost periodic solutions
of the system  ({\ref{e2n}).  The last section concerns with biasymptotic equivalence problem and
the existence of biasymptotically  almost periodic solutions of   (\ref{1}). In addition, examples are
given to illustrate the results.

\section{Asymptotic equivalence of linear systems and asymptotically almost periodic solutions}

Let $X(t)$, $X(0) = I$, be a fundamental matrix solution of  (\ref{2}).  Setting $y = X(t) u$,  we
easily see from (\ref{1}) that
\begin{eqnarray}
&& u' = P(t)u,
\label{e4}
\end{eqnarray}
where $P(t)= X^{-1}(t)B(t)X(t).$

Assume that
\begin{itemize}
\item[$\rm (C_1)$] $\displaystyle \int_{0}^{\infty} ||P(t)||dt < \infty.$
\end{itemize}

The following lemma has been obtained  R\'ab   \cite{mr1}, \cite{mr2}, for which we include a proof  for convenience.

\begin{lemma}  If  $\rm (C_1)$ is valid then the matrix differential equation
\begin{eqnarray}
&& \Psi' = P(t)(\Psi + I)
\label{e5}
\end{eqnarray}
has a solution $\Psi(t)$ which satisfies $\Psi(t) \rightarrow 0$ as $t \rightarrow \infty.$
\label{lem1}
\end{lemma}
{\it Proof}. \rm  Construct a sequence of $n \times  n$ matrices $\{\Psi_k\}$ defined on $R_+ = [0, \infty)$ as follows:
\begin{eqnarray*}\Psi_0(t) = I,\quad  \Psi_k(t) = - \int_{t}^{\infty} P(s) \Psi_{k-1}(s)\,ds\;\;
\mbox{for $k=1,2,\ldots$}.\end{eqnarray*}
Fix  $\epsilon$, $0<\epsilon<1.$ In view of ($\rm C_1$) there exists a $t_1>0$ such that
\begin{eqnarray*}\int_{t}^{\infty} ||P(s)||ds < \epsilon\quad \mbox{for all $t > t_1$.}\end{eqnarray*}
 It follows that
$||\Psi_k(t)|| < \epsilon^k, \quad k\in N,$
and consequently the series $\sum_{k=1}^{\infty} \Psi_k(t)$ is convergent uniformly for
$t \in [t_1, \infty).$ Letting  $\Psi(t) =  \sum_{k=1}^{\infty} \Psi_k(t)$,
one can easily check that $\Psi$ satisfies
\begin{eqnarray}
&& \Psi(t) = - \int_{t}^{\infty} P(s) [I+\Psi(s)]ds
\label{e6}
\end{eqnarray}
and hence it is a solution (\ref{e5}). From (\ref{e6}) it also
follows that   $\Psi \rightarrow 0$ as $t \rightarrow \infty,$
which completes that proof.

We may assume that
\begin{itemize}
\item[($\rm C_2$)] $\lim_{t \rightarrow \infty} X(t)\Psi(t) = 0.$
\end{itemize}

\begin{theorem} Suppose that conditions $(\rm C_1)$ and  $(\rm C_2)$  hold.
Then  (\ref{1}) and  (\ref{2}) are  asymptotically  equivalent.
\label{t1}
\end{theorem}
 {\bf Proof.}  Let $t$ be sufficiently large, $t\geq t_1$ say. In view of (\ref{e6}) we see that the function
$u(t) = [I + \Psi(t)]c$, $c\in \mathbb R^n$, is a solution of  (\ref{e4})
defined on $[t_1,\infty)$   and hence
\begin{equation}
y(t)=X(t)[I+\Psi(t)]c
\label{rep2} \end{equation}
 is a solution of (\ref{1}) .

Since $\Psi(t) \rightarrow 0$ as $t \rightarrow  \infty$,  there exists a
$t_2 > t_1$ such that  $I +\Psi(t_2)$ is nonsingular.
Let $x^0 = X(t_2)c$ and $y^0 = X(t_2)(I + \Psi(t_2))c$.  Denote  by   $y(t,c) = y(t, t_2, x^0)$ and
$ x(t,c) = x(t, t_2, y^0)$  the solutions
of (\ref{1}) and (\ref{2}) satisfying $x(t_2)=x^0$ and $y(t_2)=y^0$, respectively.
Now, because of the existence and uniqueness of solutions of  linear differential equations
and the fact that $I +\Psi(t_2)$ is nonsingular,
the relation $$y^0 = X(t_2)[I + \Psi(t_2)]X^{-1}(t_2)x^0$$ defines an isomorphism between
solutions $x(t)$ of (\ref{2}) and $y(t)$ of (\ref{1})
such that $$y(t) = x(t) + X(t)\Psi(t)c$$ for $t>t_1$. The last equality, in view of  $(\rm C_2)$, completes the proof.


\begin{corollary} Suppose that  the system (\ref{2}) has  a $k-$parameter  ($k \leq n$) family $\sigma_2$ of
almost periodic solutions, and that the  conditions
$\rm (C_1)$, $\rm (C_2)$ are satisfied.
Then there exists  a $k-$parameter family  $\sigma_1 $ of asymptotic almost periodic solutions of  (\ref{1}),
and $\sigma_1$  is isomorphic $\sigma_2.$
\label{t2}
\end{corollary}

\begin{example} \rm
Consider the systems
\begin{eqnarray}
&&  x'' - 2(t+1)^{-2} x = 0
\label{e21b}
\end{eqnarray}
and
\begin{eqnarray}
&&  y'' - [2(t+1)^{-2}+b(t)]y = 0,
\label{e22b}
\end{eqnarray}
where  $b(t)$ is a continuous function defined on $R_+$.
We assume that there exist real numbers $K_1> 0$ and $\alpha>0$ such that
\begin{eqnarray}
\mid b(t)\mid < K_1e^{-\alpha\, t}\quad \mbox{ for all $t\in R_+.$}\label{c1}
\end{eqnarray}
Notice that (\ref{e21b}) has solutions $x_1(t)=(t+1)^2$ and $x_2(t)=(t+1)^{-1}$.

If we transform the above second order equations into systems of the form (\ref{1}) and (\ref{2}) we identify that
$$ A = \left[
 \begin{array}{cc}
     0 & 1 \\
    2/(t+1)^{2} & 0
 \end{array}
\right]
\;\; \mbox{and $\;
 B = \left[
 \begin{array}{cc}
     0 & 0 \\
    b(t) & 0
 \end{array}
\right]$}.
$$

 It is easy to see that for a
given $\epsilon>0$ there exists $K>0$ such that
\begin{eqnarray}
&& \mid\mid P(t)\mid\mid \leq K e^{(-\alpha+\epsilon) t} \quad \mbox{for all $t\in R_+$. }
\label{e11b}
\end{eqnarray}

Fix $\epsilon$ so that $\beta:=\alpha+\epsilon <0$. Then (\ref{e4}) is satisfied, i.e.,
\begin{eqnarray*}
&& \int_{0}^{\infty}||P(t)||\,dt <\infty,
\label{e13n}
\end{eqnarray*}
and
\begin{eqnarray}
&& \mid\mid \Psi(t)\mid\mid \leq  e^{K|\beta|^{-1}e^{\beta t}} - 1.
\label{e12b}
\end{eqnarray}
Moreover, using (\ref{e11b}) and (\ref{e12b}) one can show that $X(t)\Psi(t) \rightarrow 0$ as
$ t\rightarrow \infty$, i.e.,  $\rm (C_2)$ holds.

Since the conditions of
Theorem \ref{t1} are fulfilled,  we may conclude that  (\ref{e21b})  and (\ref{e22b}) are
asymptotically equivalent whenever (\ref{c1}) holds.

\end{example}

\begin{example}\rm
Let  $b(t)$ be a continuous function such that $|b(t)| \leq K_1  e^{-\alpha t}$ for all $t\in R_+$
for some $\alpha>0$, $K_1>0$, and $C\in \mathbb R^{5\times 5}$. Consider
\begin{eqnarray}
&& y' = (A + B(t))y,
\label{e11}
\end{eqnarray}
where
\begin{displaymath}
 A= \left( \begin{array}{ccccc}
0 & 1 & 0 & 0 & 0\\
-1 & 0 & 0 & 0 & 0\\
 0 & 0 & 0 & \pi & 0 \\
0 & 0 & -\pi  & 0 & 0 \\
0 & 0 & 0 & 0 & \beta \end{array} \right),
\end{displaymath}
$B(t) = b(t)C,$  and $\beta >0$ satisfies $\alpha-2\beta>0$.

 The associated  equation $x' = Ax$ has a fundamental matrix
\begin{displaymath}
X(t) = \left( \begin{array}{ccccc}
  \cos t & \sin t &  0  &  0 &0 \\
-\sin t &  \cos t &  0   &  0 & 0 \\
 0 & 0  & \cos \pi t   & \sin \pi t & 0\\
0 & 0  & \sin \pi t  & \cos \pi t & 0 \\
 0 & 0 &  0 &  0 & e^{\beta t} \end{array} \right).
\end{displaymath}
The equality $$P(t) = X^{-1} (t) B(t) X(t) = b(t)X^{-1} (t)C X(t)$$ implies that there exists a $K>0$
such that $||P(t)|| \leq K e^{-(\alpha - \beta) t}$ for all $t\in R_+$. Therefore,
$\rm (C_1)$ is valid. We also see that
\begin{eqnarray*}||\Psi(t)|| &\leq & \sum_{k=1}^{\infty} \frac{(K e^{- (\alpha- \beta) t})^k}{ (\alpha-\beta)^k k!} =
1 - e^{K(\alpha-\beta)^{-1} e^{-(\alpha -\beta)
t}},\end{eqnarray*} and hence $X(t)\Psi(t) \rightarrow 0$ as $ t
\rightarrow \infty.$ In view of Corollary \ref{t2} we conclude
that system (\ref{e11}) has  a $4-$ parameter family of
asymptotically almost periodic solutions. More precisely they are
asymptotically "quasiperiodic" solutions and every such solution
has a torus as the $\omega-$ limit set.
\end{example}

\section{Asymptotic equivalence of linear and quasilinear  systems}

Let  $\alpha = \min_j  \Re \lambda_j$ and $\beta = \max_j  \Re \lambda_j$,
where $\Re \lambda_j$ denotes the real part of the eigenvalue $\lambda_j$ of the matrix $C$.
Let  $m_{\alpha}$  and $m_{\beta}$  be the  maximum of degrees of elementary divisors of $C$ corresponding  to
eigenvalues with real part equal to $\alpha$ and  $\beta$, respectively. Clearly, there exist
constants $\kappa_1,
\kappa_2$  such that
$$ ||e^{Ct}|| \leq \kappa_1 t^{m_{\beta}-1}e^{\beta t} \quad \mbox{and $
||e^{-Ct}|| \leq \kappa_2 t^{m_{\alpha} - 1} e^{-\alpha t}$}$$
for all $t\in R_+ = [0, \infty)$.

The following conditions are to be assumed:
\begin {itemize}
\item[$\rm (C_3)$]  $\displaystyle  || f(t,x_1)- f(t,x_2)||\leq
\eta(t)||x_1-x_2|| $ for all $ (t,x_1),(t,x_2) \in R_+\times R^n$,
and for some nonnegative function $\eta(t)$ defined on $R_+$;
\item[$\rm (C_4)$]  $ \displaystyle L:=\int_0^{\infty}
t^{m_{\beta} + m_{\alpha} - 2} e^{(\beta - \alpha)t} \eta(t)\,dt <
\infty.$
\end{itemize}
\begin{lemma}  If $\rm (C_3)$ and  $\rm (C_4)$ are valid, then every solution of
\begin{eqnarray}
&& u' = e^{-Ct} f(t, e^{Ct}u)
\label{e4n}
\end{eqnarray}
is bounded on $R_+$ and for each solution $u$ of (\ref{e4n})  there
 exists a constant vector
$c_u \in \mathbb R^n$ such that  $u(t) \rightarrow c_u$ as $t \rightarrow \infty.$
\label{lem11}
\end{lemma}
{\it Proof.} Let $u(t)=u(t,t_0,u_0)$ denote the solution of (\ref{e4n}) satisfying $u(t_0)=u_0$, $t_0\geq 0$.  It is clear that
\begin{eqnarray*}
&& u(t) = u_0 + \int_{t_0}^{t} e^{-Cs} f(s, e^{Cs}u(s))\,ds,\quad t\geq t_0.
\end{eqnarray*}
By using  $\rm (C_3)$ and $f(t,0)=0$, we see that
$$|u(t)| \leq |u_0|+k_1\int_{t_0}^{t} s^{m_{\beta} + m_{\alpha} - 2} e^{(\beta - \alpha)s} \eta(s)|u(s)|\,ds,
\quad t\geq t_0$$ for some $k_1>0$.
In view of $\rm (C_4)$  and Gronwall's inequality, we have
$$|u(t)| \leq |u_0|\,e^{\displaystyle k_1\int_{t_0}^{t} s^{m_{\beta} + m_{\alpha} - 2} e^{(\beta - \alpha)s} \eta(s)\,ds}\leq |u_0|\, e^{k_1L}<\infty,\quad
t\geq t_0.$$
Let $M_0=\max\{|u(t)|:\, t\in [0,t_0]\}$ and $M=\max \{M_0,|u_0|\, e^{k_1 L }\}$.
Then we have $|u(t)|\leq M$ for all $t\in R_+$.

To prove the second part of the theorem, we first note that
$$|\int_{t_0}^{t} e^{-Cs} f(s, e^{Cs}u(s))\,ds| \leq Mk_1\int_0^{\infty} t^{m_{\beta} +
m_{\alpha} - 2} e^{(\beta - \alpha)t} \eta(t) dt < \infty.$$  So we may define
$$c_u=u_0+\int_{t_0}^{\infty} e^{-Cs} f(s, e^{Cs}u(s))\,ds.$$
It follows that
\begin{eqnarray*}
&& u(t) = c_u - \int_{t}^{\infty} e^{-Cs} f(s, e^{Cs}u(s))ds,
\end{eqnarray*}
which completes the proof.

The following lemma can be easily justified by a direct substitution.
\begin{lemma}  If $y(t)$  is a solution of (\ref{e2n}), then there is a solution $u(t)$ of (\ref{e4n}) such that
\begin{eqnarray}
&& y(t) =  e^{Ct}u(t).
\label{e8}
\end{eqnarray}
Conversely, if $u(t)$ is a solution of  (\ref{e4n})
then $y(t)$ in (\ref{e8}) is a solution of (\ref{e2n}).
\label{lem33}
\end{lemma}

\begin{theorem}
If conditions $\rm (C_3)$ and $\rm (C_4)$ are satisfied, then
every solution $y(t)$ of (\ref{e2n}) possesses an asymptotic
representation of the form
\begin{eqnarray*}
&& y(t) =  e^{C\,t} [ c + o(1)],
\end{eqnarray*}
where $c \in \mathbb R^n$ is a constant vector and for a solution $u(t)$ of
(\ref{e4n}), $$o(1) =  -\int_{t}^{\infty} e^{-Cs} f(s, e^{Cs}u(s))\,ds.$$
\label{t1n}
\end{theorem}

{\it Proof.} The proof follows from Lemma \ref{lem11} and Lemma \ref{lem33}.

\begin{theorem}
Assume that $\rm (C_3)$ and $(C_4)$ are fulfilled, and
\begin{itemize}
\item[$\rm (C_5)$] $ \displaystyle \lim_{t\to\infty}\int_{t}^{\infty} (s-t)^ {m_{\alpha} - 1}s^{m_{\beta} -1}
  e^{\alpha(t-s)}     e^{\beta s} \eta(s) ds =0.$
\end{itemize}
Then (\ref{e2n}) and  (\ref{e1n}) are asymptotically equivalent.
\label{ya}
\end{theorem}

{\it Proof.}  In view of Lemma \ref{lem11} we see that
\begin{eqnarray*}
y(t) &=&   e^{Ct}[c_u - \int_{t}^{\infty} e^{-Cs} f(s, e^{Cs}u(s))ds]  \nonumber\\
&=&x(t) - \int_{t}^{\infty} e^{C(t-s)} f(s, e^{Cs}u(s))ds,
\end{eqnarray*}
where  $x(t) = e^{Ct}c_u $ is a solution of  (\ref{e1n}) and $u(t)
= u(t,t_0,u_0)$ is a solution of  (\ref{e4n}). It is clear that a
given $u_0$ results in a one-to-one correspondence between $x(t)$
and $y(t).$ In view of  $\rm (C_5)$, we also see that $x(t) - y(t)
\rightarrow 0 $ as $t \rightarrow \infty,$ which completes the
proof of the theorem.

In \cite{ya}, Yakubovich proved that if
\begin{eqnarray}
&&\lim_{t\to\infty} \int_t^{\infty}s^{m_{\beta}+m_{\alpha} - 2}
 e^{\beta s} \eta(s) ds = 0
\label{yak}
\end{eqnarray}
then  (\ref{e2n}) and  (\ref{e1n}) are asymptotically equivalent.
It is clear that if $\alpha>0$ then  condition $\rm (C_5)$  is weaker than (\ref{yak}).

The following assertion  is a simple corollary of the Theorem \ref{ya}

\begin{corollary} Suppose that conditions $\rm (C_3)$, $(\rm C_4)$,  $\rm (C_5)$ hold, and that
system (\ref{e1n}) has a $k-$ parameter ($k \leq n$) family $\gamma_1$ of
almost periodic solutions.  Then  (\ref{e2n}) admits a  $k-$parameter  family $\gamma_2$ of
asymptotically almost periodic solutions, and  $\gamma_1$  is homeomorphic $\gamma_2.$
\label{0}
\end{corollary}

\section{Biasymptotic equivalence of linear systems and biasymptotically almost periodic solutions}

With regard to systems (\ref{1}) and (\ref{2}) we shall make use of the following conditions:
\begin{itemize}
\item[$(\rm C_6)$] $A(-t) = - A(t)$ for all $t \in \mathbb R.$
\item[$(\rm C_7)$] $B(-t) = B(t)$ for all $t \in \mathbb R.$
\end{itemize}

We will rely on the following two lemmas. The first lemma is
almost trivial.

\begin{lemma} If  $\rm (C_6)$ is satisfied then $ X(-t) = X(t)$ for all $t \in \mathbb R,$ and if in addition
$\rm (C_7)$ holds then   $ P(-t) = -P(t)$ for all $t \in \mathbb R$.
\label{lem1a}
\end{lemma}


\begin{lemma}  Assume that conditions $\rm (C_1)$, $\rm (C_6)$, $\rm (C_7)$ are valid.
Then (\ref{e5}) has a solution $\Psi(t)$ which satisfies
$\Psi(-t) = \Psi(t)$  for all  $t \in \mathbb R$ and $\Psi \rightarrow 0$ as $t \rightarrow \infty.$
\label{lem4}
\end{lemma}

{\it Proof.} \rm By Lemma   \ref{lem1} there exists a solution $\Psi_+(t)$ of (\ref{e5}) which is defined
for $t \geq t_1$ and satisfies $\Psi_+ \rightarrow 0$ as $t \rightarrow \infty.$

Using $P(-t)=-P(t)$ we see that
$$\int_{-\infty}^{-t_1} ||P(s)||\,ds  = \int_{t_1} ^{\infty}||P(s)||\,ds< \epsilon.$$
We may define a sequence of $n \times  n$ matrices $\{\tilde \Psi_k\}$ for $t\in (-\infty,0]$ as follows:
\begin{eqnarray*} \tilde{\Psi}_0(t) = I,\quad   \tilde{\Psi}_k =  \int_{-\infty}^{t} P(s)\tilde{\Psi}_{k-1} (s)\,ds \;\; \mbox{for $ k=1,2,\ldots$}.\end{eqnarray*}
As in the proof of Lemma \ref{lem1}, the matrix function
$\Psi_-(t) = \sum_{k=1}^{\infty}  \tilde{\Psi}_k (t)$ satisfies
$$\Psi(t) =   \int_{-\infty}^{t} P(s) (I + \Psi(s)) ds$$ and hence  becomes a solution of (\ref{e5}) for $t \leq - t_1.$

On the other hand, since $\Psi_0(-t)= \tilde{\Psi}_0(t) = I$ we have
\begin{eqnarray*}\Psi_1(-t) &=& - \int_{-t}^{\infty} P(s)\Psi_0 (s) ds = \int_{t}^{-\infty} P(-s)\Psi_0 (-s) ds \\
&=&  \int_{-\infty}^{t} P(s)\tilde \Psi_0 (s) ds = \tilde \Psi_1 (t). \end{eqnarray*}
It follows by induction that $\Psi_k(-t) = \tilde \Psi_k(t)$  for all $k=0,1,2,\ldots$ and for all $t\leq - t_1$.
Hence
$\Psi_+(-t) = \Psi_-(t)$ for all $t \leq -  t_1.$
Continuing $\Psi_+$ and  $\Psi_-$ as solutions of   (\ref{e5}),  one can obtain that
$\Psi_+(-t) = \Psi_-(t)$ for all  $t\leq 0.$

Define
\begin{eqnarray*}
&&  \Psi(t) = \left\{\begin{array}{rr} \Psi_+ & \mbox{if  $t \geq 0$},\\
\Psi_- & \mbox{if  $t<0$}.
\end{array}\right.
\end{eqnarray*}
Clearly $\Psi(t)$ is a solution of (\ref{e5}) satisfying $\Psi(-t) = \Psi(t)$  for all  $t \in \mathbb R$ and $\Psi \rightarrow 0$ as $t \rightarrow \infty.$ This completes the proof.

The following results are analogous to Theorem \ref{t1} and Corollary \ref{t2}.
\begin{theorem} Suppose that $\rm (C_1)$, $\rm (C_2)$, $\rm (C_6)$, $\rm (C_7)$ are valid.
Then  (\ref{1}) and (\ref{2}) are asymptotically biequivalent.
\label{t55}
\end{theorem}

\begin{corollary} Suppose that  $\rm (C_1)$, $\rm (C_2)$, $\rm (C_6)$, $\rm (C_7)$ are valid, and that
(\ref{2}) has a $k-$ parameter ($k \leq n$) family $\nu_1$ of
almost periodic solutions.
Then  (\ref{1}) admits a  $k-$parameter  family $\nu_2$ of
biasymptotically almost periodic solutions, and   $\nu_1$  is isomorphic  $\nu_2.$
\label{t56}
\end{corollary}

\begin{example} \rm Consider the system
\begin{eqnarray}
&& y' = (A(t) + B(t))y,
\label{e12}
\end{eqnarray}
where
\begin{displaymath}
A(t) = \left( \begin{array}{cc}
\sin \pi t   & 0 \\
0 & \sin{\sqrt{5}t} \end{array} \right),
\end{displaymath}
 $B(t) = \cos t\, e^{-\alpha |t|}\,C$ with $\alpha >0$ a real number and  $C\in \mathbb R^{2\times 2}$.

Applying Corollary \ref{t56} one can conclude that every solution of  (\ref{e12})
is biasymptotically quasiperiodic.
\end{example}

\end{document}